\documentclass[a4paper,UKenglish,cleveref, autoref]{lipics-v2021}

\pdfoutput=1 
\hideLIPIcs  

 \title{Topological Inference for degree-\v{C}ech complexes}

\titlerunning{}

\author{Michael Kerber}{Institute of Geometry, Graz University of Technology, Austria}{kerber@tugraz.at}{https://orcid.org/0000-0002-8030-9299}{Austrian Science Fund (FWF) grant 10.55776/DOC183}

\author{Bria Razanaparany}{Institute of Geometry, Graz University of Technology, Austria}{bria.razanaparany@tugraz.at}{https://orcid.org/0009-0000-0880-1660}{Austrian Science Fund (FWF) grant 10.55776/DOC183}

\authorrunning{Kerber and Razanaparany}


\category{} 

\relatedversion{} 

\acknowledgements{}
\nolinenumbers 

\Copyright{Kerber, Razanaparany}

\EventEditors{John Q. Open and Joan R. Access}
\EventNoEds{2}
\EventLongTitle{42nd Conference on Very Important Topics (CVIT 2016)}
\EventShortTitle{}
\EventAcronym{CVIT}
\EventYear{2016}
\EventDate{December 24--27, 2016}
\EventLocation{Little Whinging, United Kingdom}
\EventLogo{}
\SeriesVolume{42}
\ArticleNo{23}

\usepackage{thm-restate}

 \usepackage{mathtools}
\usepackage[table]{xcolor}
 \AddToHook{cmd/appendix/before}{\crefalias{section}{appendix}}

\newcommand{\ignore}[1]{}
\newcommand{\R}{\mathbb{R}}

\DeclareMathOperator{\Vol}{Vol}

\newcommand{\Eout}{E_{\mathrm{out}}}
\newcommand{\Ein}{E_{\mathrm{in}}}

\newcommand{\eps}{\varepsilon}

\newcommand{\Cech}{\v{C}ech\ }

\DeclarePairedDelimiterX\Set[1]\{\}{%
  
  #1
}

\usepackage{tikz-cd}
\usetikzlibrary{arrows}
\usetikzlibrary{arrows.meta,positioning}

\usepackage[ruled,vlined, linesnumbered]{algorithm2e}

\begin{document}

\maketitle

\begin{abstract}
We study topological inference for the degree-\Cech bifiltration from noisy samples of a manifold.
Under a model combining uniform sampling on the manifold with uniform noise, we show that suitable choices of the radius and degree parameters recover the homotopy type of the underlying manifold with high probability.
The analysis explains how the degree parameter suppresses sufficiently distant outliers while the radius scale captures topological signals.
Our results provide theoretical support for the common intuition that degree bifiltrations carry meaningful topological information in the presence of noise.
\end{abstract}

\section{Introduction}
\emph{Topological inference}~\cite{bcy-geometric} is the task
to deduce topological properties of an unknown space from a possibly noisy
point sample that approximates that space.
The motivation behind the question is that topology offers qualitative
and robust information about the underlying space which, for instance,
facilitates the analysis of geometric properties of the space as well~--
see also the related area of \emph{manifold learning}~\cite{mz-manifold}.

Most theorems in topological inference follow the pattern popularized
by Niyogi, Smale, and Weinberger~\cite{nsw-finding}: given a ``good sampling'',
of a geometric shape
specified by a set of \emph{sampling conditions} and a suitable choice
of a \emph{scale parameter}, a simplicial approximation based on the point
sample (e.g., a \Cech- or Vietoris-Rips complex) is guaranteed to yield
the same topological invariant (e.g., same homotopy type or
isomorphic homology groups) as the underlying shape. There is a substantial
body of literature following this template~\cite{akfflsw-tight,k-j-c-w-homot-recons,c-l-manifol-reconst,h-vietoris-rips,c-o-towards-pers-recons}.
A slight variant is \emph{probabilistic inference}
where the point sample is drawn from a specified distribution and the
invariant of the underlying space is recovered with good probability~\cite{nsw-finding,n-s-w-top-learn,c-s-q-geom-infer-prob-meas}.

Most topological inference results can cope with noise in the point sample.
This usually means that the sample is not required to be directly on the shape,
but geometrically close to it (e.g., with bounded Hausdorff distance).
However, these methods usually fail in the presence of uniform noise
in the sample. Indeed, such type of noise tends to dilute the homotopy
or homology signal of the point sample more strongly and makes inference
much harder.
A natural remedy is to filter out low-density points via some outlier detection
and apply the inference method only on the remaining points.

The presence of uniform noise is also a main motivation for
\emph{multi-parameter persistent homology}~\cite{bl-introduction}.
Several extensions of \Cech- and Vietoris-Rips complexes were suggested
in which a second parameter controls which data points are treated as outliers
and ignored in the reconstruction. A simple example is the \emph{degree-\Cech complex} that will be studied in this work. The general motto is
that the topological signal of the underlying manifold should persist over a
large range of parameters, making a multi-scale analysis of the $2$-parameter
space for the point cloud meaningful.

\subparagraph{Contribution.}
We prove probabilistic topological inference theorems
on degree-\Cech complexes
as well as on core bifiltrations~\cite{bbgs-core} as an intermediate step.
For that, we define a simple mixed distribution
consisting of the signal, a uniform point sample of a manifold embedded in the $d$-dimensional unit cube, and the noise,
a uniform point sample on that unit cube. We then analyze for which choices
of parameters the degree-\Cech complex reconstructs the homotopy type
of the manifold (with fixed probability $1-\delta$).
Our approach is to ensure that all points sufficiently far away
from the manifold are marked as outliers with our choice of parameters,
and the remaining points recover the homotopy type of the manifold,
using an existing inference result by Attali et al.~\cite{akfflsw-tight}.
While the main result is technical and depends on various parameters such as
the reach and dimension of the manifold, the number of sample points,
the signal-to-noise ratio of the distribution and other things,
all major ideas of our proof can be explained with minimal technical
exposition, and we will start with such an exposition in Section~\ref{sec:overview}.

We see our result as a step towards extending topological inference to the
case of bi-filtered complexes which seems natural and useful for 
application scenarios. At the same time, our work supports
the folklore claim that the bifiltrations
considered in multi-parameter persistence carry a meaningful topological
signal.

\subparagraph{Related work.}
The work most closely related to ours is by Balakrishnan et al~\cite{balakrishnan-minimax}:
they consider a noise model very similar to ours
(called ``clutter noise'' in their work)
and show that when removing outlier vertices, the \Cech complex for
a specific radius yields the same homology as the sampled shape. Their
approach to arrive at this result is in large part similar to our exposition
in Section~\ref{sec:overview} but differs in some details; we postpone
a more precise discussion of the differences to the end of Section~\ref{sec:overview}.

Another related work is Rolle's study of the degree-Rips bifiltration for an annulus with outliers~\cite{a-role-deg-rips-annulus}. 
In that setting, the degree-Rips complexes are computed explicitly up to homotopy type.
In particular, this work advocates a homological inference perspective for understanding why degree bifiltrations can recover the topological structure in the presence of outliers.
Our work follows the same perspective but addresses a more general setting.
Much work has been done towards bifiltration, that has primarily focused on stability and robustness.
Blumberg and Lesnick established stability results for several bifiltrations that allow the addition and removal of outliers~\cite{bl-l-stability}.
Degree bifiltration constructions are also related to density based clustering methods such as DBSCAN and its variants~\cite{r-s-stable-clustering}.

\section{Overview of our approach.}
\label{sec:overview}

We describe the major ideas of how we obtain our main result in this section.
We will assume some basic knowledge in combinatorial topology and introduce
formal concepts only briefly and to an extent that is necessary to state
the major steps informally.

\begin{figure}
  \centering
  \includegraphics[width=10cm]{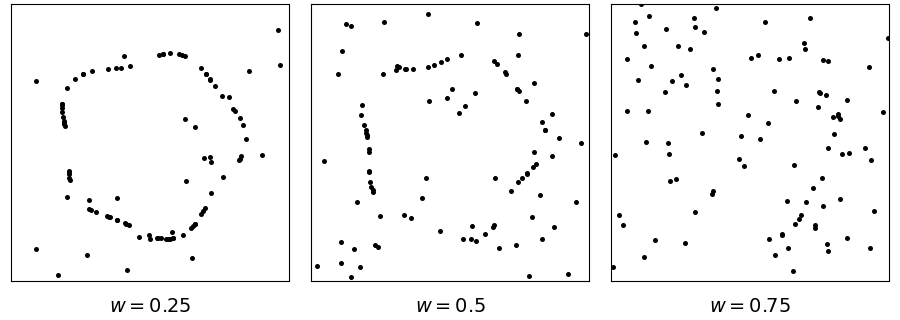}
  \caption{Sample of 100 points with signal-to-noise-ratio of $w=0.25$, $w=0.5$ and $w=0.75$.}
  \label{fig:samples}
\end{figure}

\subparagraph{The point sample.}
We fix $M$ as a smooth compact manifold without boundary, embedded into the unit
cube $[0,1]^d$. We fix $n$ as the number of points in our sample and $w\in[0,1)$
as the \emph{signal-to-noise-ratio} and generate the point sample $P$
by drawing $\lceil (1-w)n \rceil$ points uniformly at random
from $M$ (without noise)
and $\lfloor wn\rfloor$ points uniformly at random from $[0,1]^d$
(we ignore the rounding to integers for simplicity in the analysis).
See Figure~\ref{fig:samples} for an illustration. In what follows,
we refer to $P$ as a $w$-sample of $M$ for brevity.

\begin{figure}[thbp]
    \centering
\begin{tikzpicture}[scale=1]
  \coordinate (A) at (6.11, 8.45);
  \coordinate (C) at (6.11, 8.45);
  \coordinate (D) at (5.10, 8.41);
  \coordinate (E) at (5.10, 8.41);
  \coordinate (F) at (4.63, 7.62);
  \coordinate (G) at (4.63, 7.62);
  \coordinate (H) at (6.55, 8.06);
  \coordinate (I) at (6.55, 8.06);N
  \coordinate (J) at (5.50, 8.50);
  \coordinate (K) at (5.50, 8.50);
  \coordinate (L) at (5.03, 6.77);
  \coordinate (M) at (5.03, 6.77);
  \coordinate (P) at (6.98, 6.96);
  \coordinate (Q) at (6.98, 6.96);
  \coordinate (R) at (6.33, 6.36);
  \coordinate (S) at (6.33, 6.36);
  \coordinate (T) at (6.44, 6.87);
  \coordinate (U) at (6.44, 6.87);
  \coordinate (V) at (6.75, 6.31);
  \coordinate (W) at (6.75, 6.31);
  \coordinate (B) at (5.00, 8.01);
  \coordinate (N) at (5.00, 8.01);
  \draw[white, ultra thin, fill=green, fill opacity=0.3] (A) circle (0.56);
  \draw[white, ultra thin, fill=green, fill opacity=0.3] (D) circle (0.56);
  \draw[white, ultra thin, fill=green, fill opacity=0.3] (F) circle (0.56);
  \draw[white, ultra thin, fill=green, fill opacity=0.3] (H) circle (0.56);
  \draw[white, ultra thin, fill=green, fill opacity=0.3] (J) circle (0.56);
  \draw[white, ultra thin, fill=green, fill opacity=0.3] (L) circle (0.56);
  \draw[white, ultra thin, fill=green, fill opacity=0.3] (P) circle (0.56);
  \draw[white, ultra thin, fill=green, fill opacity=0.3] (R) circle (0.56);
  \draw[white, ultra thin, fill=green, fill opacity=0.3] (T) circle (0.56);
  \draw[white, ultra thin, fill=green, fill opacity=0.3] (V) circle (0.56);
  \filldraw[very thin, fill=blue, fill opacity=0.3] (R) -- (V) -- (P) -- cycle;
  \filldraw[very thin, fill=blue, fill opacity=0.3] (R) -- (T) -- (P) -- cycle;
  \draw[very thin] (T) -- (V);
  \filldraw[very thin, fill=blue, fill opacity=0.3] (F) -- (B) -- (D) -- cycle;
  \filldraw[very thin, fill=blue, fill opacity=0.3] (B) -- (J) -- (D) -- cycle;
  \draw[very thin] (L) -- (F);
  \draw[very thin] (J) -- (A);
  \draw[very thin] (A) -- (H);
  \draw[white, ultra thin, fill=green, fill opacity=0.3] (B) circle (0.56);
  \fill[black] (A) circle (1.5pt);
  \fill[black] (D) circle (1.5pt);
  \fill[black] (F) circle (1.5pt);
  \fill[black] (H) circle (1.5pt);
  \fill[black] (J) circle (1.5pt);
  \fill[black] (L) circle (1.5pt);
  \fill[black] (P) circle (1.5pt);
  \fill[black] (R) circle (1.5pt);
  \fill[black] (T) circle (1.5pt);
  \fill[black] (V) circle (1.5pt);
  \fill[black] (B) circle (1.5pt);
\end{tikzpicture}
    \caption{Point set $X$ and its corresponding \Cech complex at scale radius $r$.}
    \label{fig:cech_illust}
\end{figure}
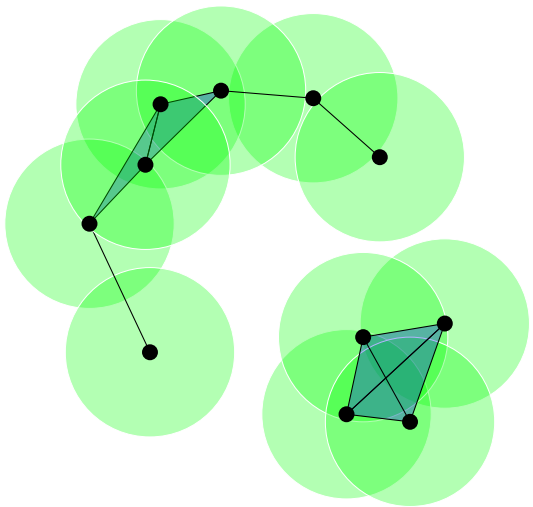

\subparagraph{Degree-\Cech complexes.}
Recall that the \Cech complex $C_r(X)$ of a point set $X$ is the nerve
of the balls of radius $r$ centered at the points of $X$. In other words,
a subset of points $\{x_0,\ldots,x_k\}$ in $X$ forms a simplex in $C_r(X)$
if and only if the $r$-balls around $x_0,\ldots,x_k$ have a common point of intersection (Figure~\ref{fig:cech_illust}).
The Nerve Theorem~\cite{bkrr-unified-nerve} asserts that $C_r(X)$ is homotopically equivalent to the union of $r$-balls
around $X$ which we also denote as $X^r$ in the paper.

We define the \emph{degree-\Cech complex} next. In fact, we will define a slight
generalization of it which simplifies the exposition. Fix a parameter
$c\in(0,1]$. For a point $x\in X$ and $r\geq 0$, let $N_{cr}(x)$ denote
the number of points in $X$ within distance $2cr$ from $x$; we call
such points the \emph{neighbors of $x$}, with the understanding
that $x$ is always neighbor to itself. Note that if $c=1$, the
number of neighbors equals the degree of $x$ in the \Cech complex.

The \emph{$c$-degree-\Cech complex} $C_{k,r}^{(c)}$ for
radius $r$ and an integer $k\geq 1$ is the \Cech complex $C_r(X_{k,cr})$,
where $X_{k,cr}\subseteq X$ is the subset of points that have at least $k$
neighbors. Note that for $k=1$, $X_{k,cr}=X$, and we retain the ordinary
\Cech complex. Also, the usual degree-\Cech complex is obtained
by the choice $c=1$.
We point out that the $c$-degree-\Cech complex is a special
case of the recently introduced \emph{core bifiltration}~\cite{bbgs-core};
we still prefer the term ``$c$-degree-\Cech complex'' for clarity of notation.

\subparagraph{Our approach for $c<\frac{1}{2}$.}
Given a $w$-sample of $M\subset\R^d$ and some $c<\frac{1}{2}$, we want to argue that the
$c$-degree-\Cech complex will capture the homotopy type of $M$ for a range
of choices of $r$ and $k$.
There are two major ingredients to obtain this result.
For that,
we call the sample points obtained from the uniform sampling of $M$ the
\emph{inliers} and the other points the \emph{outliers}.
The first ingredient is a homotopical inference result for \Cech complexes
by Attali et al.~\cite{akfflsw-tight}. We show that, sampling sufficiently
many points, the \Cech complex of the inliers satisfies
the conditions of their theorem and yields the right homotopy
type for a certain range of radii. Moreover, their inference result is
stable enough so that the homotopy type is still preserved when adding
any outliers to the sample, as long as their distance to the manifold
is smaller than $2cr$ (for any $c<\frac{1}{2}$).

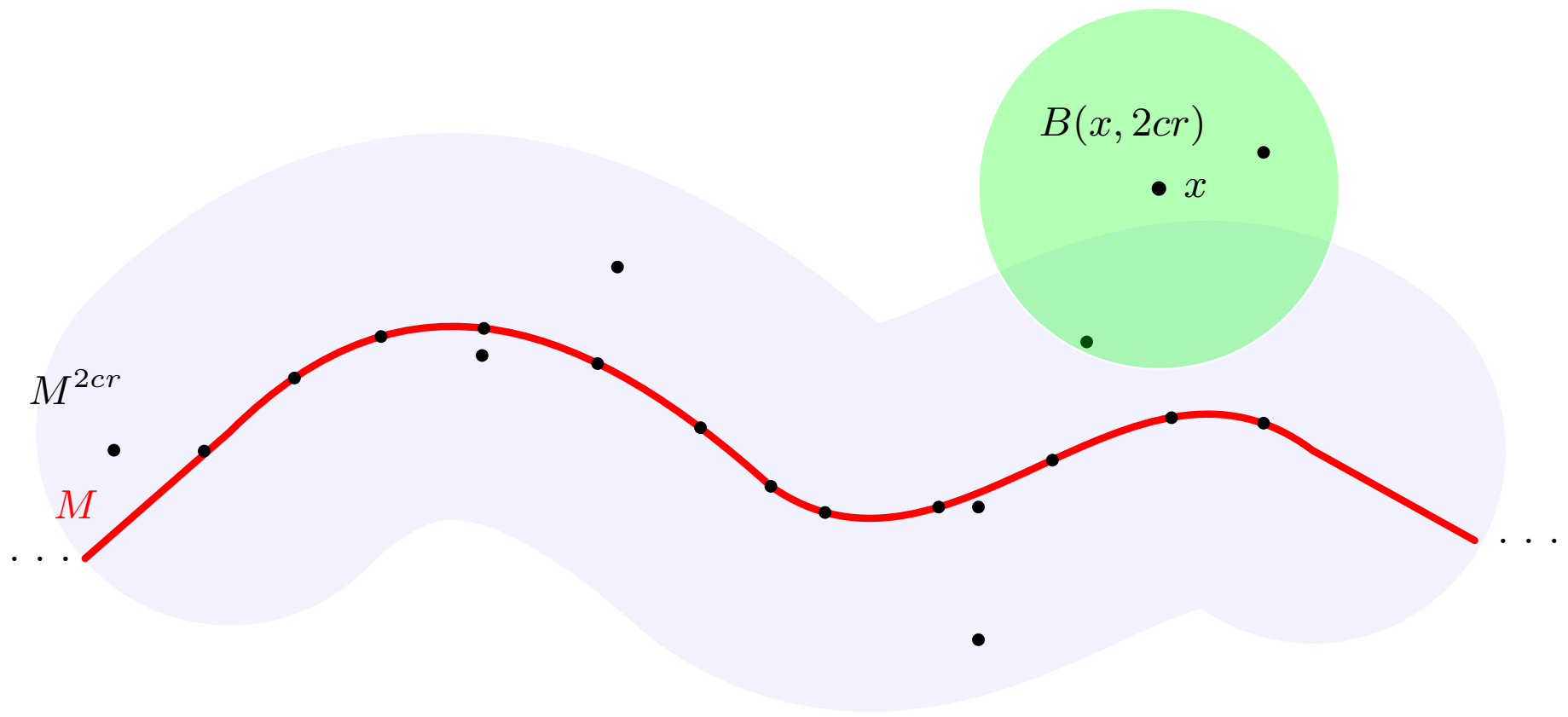
\begin{figure}[thbp]
\centering
\begin{tikzpicture}[scale=0.7, >=Latex,every node/.style={font=\small}]

\draw[
    blue!5,
    line width=3cm,
    line cap=round,
    line join=round
]
(-6,0)
.. controls (-4,2.0) and (-2,1.2) .. (0,-0.6)
.. controls (2,-2.0) and (4,1.3) .. (6,-0.2);

\draw[
    red,
    ultra thick,
    line cap=round
]
(-7.6,-1.4)
.. controls (-6,0) and (-6,0) .. (-6,0)
.. controls (-4,2.0) and (-2,1.2) .. (0,-0.6)
.. controls (2,-2.0) and (4,1.3) .. (6,-0.2)
.. controls (7.8,-1.2) .. (7.8,-1.2);

\node[black] at (8.4,-1.2) {$.\ .\ .$};
\node[black] at (-8.1,-1.4) {$.\ .\ .$};

\node[red] at (-7.7,-0.8) {$M$};
\node[black] at (-7.7,0.5) {$M^{2c r}$};

\fill[black] (-7.28, -0.2) circle (2pt);
\fill[black] (-6.28, -0.21) circle (2pt);

\fill[black] (-5.28, 0.60) circle (2pt);
\fill[black] (-4.32, 1.06) circle (2pt);
\fill[black] (-3.18, 1.15) circle (2pt);
\fill[black] (-1.92, 0.76) circle (2pt);
\fill[black] (-0.78, 0.05) circle (2pt);

\fill[black] (0.60,-0.89) circle (2pt);
\fill[black] (1.86,-0.83) circle (2pt);
\fill[black] (3.12,-0.31) circle (2pt);
\fill[black] (4.44, 0.16) circle (2pt);
\fill[black] (5.46, 0.10) circle (2pt);

\coordinate (xin) at (-3.2,0.85);
\fill[black] (xin) circle (2pt);
\node[above left=2pt of xin]{};
\coordinate (z) at (0,-0.6);
\coordinate (xnear) at (3.5,1.0);
\fill[black] (z) circle (2pt);
\fill[black] (xnear) circle (2pt);
\node[right=2pt of xnear]{};
\node[below right=1pt of z]{};
\node[right] at (1.9,1.8){};
\node[right] at (1.9,-0.1){};
\coordinate (xfar) at (4.3,2.7);
\draw[white, thin,  fill=green, fill opacity=0.3](xfar) circle (2cm);

\filldraw[black] (xfar) circle (2pt);

\fill[black] (2.3,-0.83) circle (2pt);
\fill[black] (2.3,-2.3) circle (2pt);
\fill[black] (-1.7,1.83) circle (2pt);
\fill[black] (5.46, 3.10) circle (2pt);


\coordinate (xproj) at (4.4, 0.05);


\node[right=2pt of xfar]{$x$};
\node[above] at (3.9,3.0) {$B(x,2c r)$};
\node at (4.3,1.35){};
\end{tikzpicture}
\caption{
Illustration of points near and far from the manifold $M$.
The ball centered at a far outlier point with radius $2cr$ does not intersect with the manifold $M$. 
A far outlier point can have a near outlier point as a neighbor.
}
\label{fig:neigh_out_in_illust}
\end{figure}

In the light of the last statement, we partition the outliers into
\emph{near outliers} that have distance smaller than $2cr$ to $M$,
and \emph{far outliers} with larger distance to $M$. See Figure~\ref{fig:neigh_out_in_illust}.

Recalling that $C^{(c)}_{k,r}$ is the \Cech complex of radius $r$
on the point set with at least $k$ neighbors, we only have to make sure that
(i) all inliers in the sample have at least $k$ neighbors (to ensure
homotopy equivalence) and (ii) no far outlier has $k$ or more neighbors
(because the inference result would not apply anymore if such a point is
present). Again, the near outliers might or might not have enough neighbors
to be included, but their presence or absence
does not influence the reconstruction.

The remaining challenge is to bound the number of neighbors of inliers
and far outliers to satisfy (i) and (ii) above. While technically
a bit involved, the approach is standard. Note that 
the $c$-degree-\Cech complex only counts neighbors that are in distance $2cr$
and hence, no far outlier has an inlier as neighbor (here, the assumption $c<\frac{1}{2}$
is crucial). Hence, counting its neighbors is a simple analysis of a Bernoulli
distribution, and we exploit that the number of neighbors is sharply
concentrated around the mean. A similar argument
applies for the inliers: As we need a lower bound, it suffices to count
the number of inlier neighbors, and it is enough to count that number
in a $2cr$-ball on $M$, as $M$ is embedded, and so that ball is contained
in the ambient Euclidean $2cr$-ball.

Following that strategy, we obtain our first main result:

\begin{restatable}{restatablethm}{MainTheoremOne}\label{thm:prob_inference_bifiltration}

Let $M$ be a $m$-dimensional manifold without boundary with positive reach $\tau$, embedded in the unit cube $[0,1]^d$. 
Let $P$ be a $w$-sample of $M$. Let $c \in (0, \frac{1}{2})$.
Fix $\delta \in (0,1)$ and suppose that 
$
n \geq \frac{8.016^m\Vol(M)}{(1-w)v_m(1-2c)^{2m}\tau^m }\left( \log\left(\frac{16.008^m\Vol(M)}{v_m(1-2c)^{2m}\tau^m }\right) + \log \frac{2}{\delta}\right)
$.
If 
\begin{align*} r \in\left[\frac{1-2c}{3}\tau, \frac{5(1-2c)}{6}\tau\right]
\end{align*}
and $k \in \mathbb{N}$ satisfies 
\begin{align*}
k \in \left (1+wnv_d(2cr)^d + \sqrt{\frac{wn}{2}\log\frac{2n}{\delta}}, \frac{(1-w)nv_m}{\Vol(M)}\left(\frac{47cr}{24}\right)^m  - \sqrt{\frac{(1-w)n}{2}\log\frac{2n}{\delta}}\right]
\end{align*}
 then, the $c$-degree-\Cech bifiltration $C_{k,r}^{(c)}(P)$ has the homotopy type of $M$, with probability at least $1-\delta$.
\end{restatable}

In the theorem, $v_d$ is the volume of the $d$-dimensional unit ball.
The dependence on the reach of the manifold, and the admissible range of $r$ comes from the inference
theorem that we employ from~\cite{akfflsw-tight}. For readers not familiar
with that notion of the reach, we discuss this in more detail in Section~\ref{sec:inference}.
The admissible range of $k$ comes from separating the degrees of inliers from those of far outliers.

\subparagraph{Extension to $c\geq \frac{1}{2}$.}
For the ``real'' degree-\Cech complex (i.e. setting $c=1$),
we follow a similar strategy:
we make sure that the inliers get enough neighbors to be
included in the degree-\Cech complex, and the analysis is the same.
For the outliers, however, a refinement is needed: if we defined
the threshold between near and far outliers to be $2cr=2r$, the inference
theorem of Attali et al.~does not apply anymore. Hence, we define the threshold to be
$\gamma r$ with some $\gamma<1$ (in fact, we fix $\gamma=\frac{49}{50}$ for simplicity).
The problem is that doing so, a far outlier $x$
with distance to the manifold between $\gamma r$ and $2r$ can have inlier
neighbors now, and counting the neighbors of such outliers is more involved.

The number of outlier neighbors of $x$ can be bounded as before.
The main idea is now: Because $x$ is ``far away'' from $M$, and because
the radius $r$ is well below the reach of the manifold,
all inlier neighbors of $x$ are inside a ball of radius strictly smaller than $2r$.
That means that $x$ will acquire fewer inlier neighbors in expectation
than an inlier point itself.
With some additional technical estimates in place, we arrive at our second main theorem.

\begin{figure}[ht]
\centering
\begin{tikzpicture}[scale=0.7, >=Latex,every node/.style={font=\small}]

\draw[
    blue!5,
    line width=3cm,
    line cap=round,
    line join=round
]
(-6,0)
.. controls (-4,2.0) and (-2,1.2) .. (0,-0.6)
.. controls (2,-2.0) and (4,1.3) .. (6,-0.2);

\draw[
    red,
    ultra thick,
    line cap=round
]
(-7.6,-1.4)
.. controls (-6,0) and (-6,0) .. (-6,0)
.. controls (-4,2.0) and (-2,1.2) .. (0,-0.6)
.. controls (2,-2.0) and (4,1.3) .. (6,-0.2)
.. controls (7.8,-1.2) .. (7.8,-1.2);

\node[black] at (8.4,-1.2) {$.\ .\ .$};
\node[black] at (-8.1,-1.4) {$.\ .\ .$};

\node[red] at (-7.7,-0.8) {$M$};
\node[black] at (-7.7,0.5) {$M^{\gamma r}$};

\fill[black] (-7.28, -0.2) circle (2pt);
\fill[black] (-6.28, -0.21) circle (2pt);

\fill[black] (-5.28, 0.60) circle (2pt);
\fill[black] (-4.32, 1.06) circle (2pt);
\fill[black] (-3.18, 1.15) circle (2pt);
\fill[black] (-1.92, 0.76) circle (2pt);
\fill[black] (-0.78, 0.05) circle (2pt);

\fill[black] (0.60,-0.89) circle (2pt);
\fill[black] (1.86,-0.83) circle (2pt);
\fill[black] (3.12,-0.31) circle (2pt);
\fill[black] (4.44, 0.16) circle (2pt);
\fill[black] (5.46, 0.10) circle (2pt);

\coordinate (xin) at (-3.2,0.85);
\fill[black] (xin) circle (2pt);
\node[above left=2pt of xin]{};
\coordinate (z) at (0,-0.6);
\coordinate (xnear) at (3.5,1.0);
\fill[black] (z) circle (2pt);
\fill[black] (xnear) circle (2pt);
\node[right=2pt of xnear]{};
\node[below right=1pt of z]{};
\node[right] at (1.9,1.8){};
\node[right] at (1.9,-0.1){};
\coordinate (xfar) at (4.3,2.7);
\draw[white, thin,  fill=green, fill opacity=0.3](xfar) circle (3.5cm);
\draw[white, thin,  fill=yellow, fill opacity=0.5](xfar) circle (1.5cm);

\filldraw[black] (xfar) circle (2pt);

\fill[black] (2.3,-0.83) circle (2pt);
\fill[black] (2.3,-2.3) circle (2pt);
\fill[black] (-1.7,1.83) circle (2pt);
\fill[black] (5.46, 3.10) circle (2pt);

\coordinate (xxin) at (-4.32, 1.06);
\draw[white, thin,  fill=green, fill opacity=0.3](xxin) circle (3.5cm);
\node[above=2pt of xxin]{$y$};
\node[above] at (-4.32,3.0) {$B(y, 2r)$};

\coordinate (xproj) at (4.4, 0.05);
\draw[
    blue,
    dashed,
    thick
]
(xproj) circle (2.15cm);
\draw[
    dashed,
    thick
]
(xproj) -- (xfar)
node[midway,right=2pt] {};


\node[right=2pt of xfar]{$x$};
\node[right] at (5.9,3.4) {$B(x,2r)$};
\node[above] at (3.9,3.0) {$B(x,\gamma r)$};
\node at (4.3,1.35){};
\end{tikzpicture}
\caption{
Illustration of points far from the manifold $M$. A point $x$ located within distance farther than $\gamma r$ from $M$ but strictly less than $2r$, may see some neighbors in the manifold within distance $2r$. The dashed blue circle represents a ball centered at the nearest projection point of $x$ in the manifold $M$ with radius  smaller than $2r$.
}
\label{fig:hard_neigh_illust}
\end{figure}
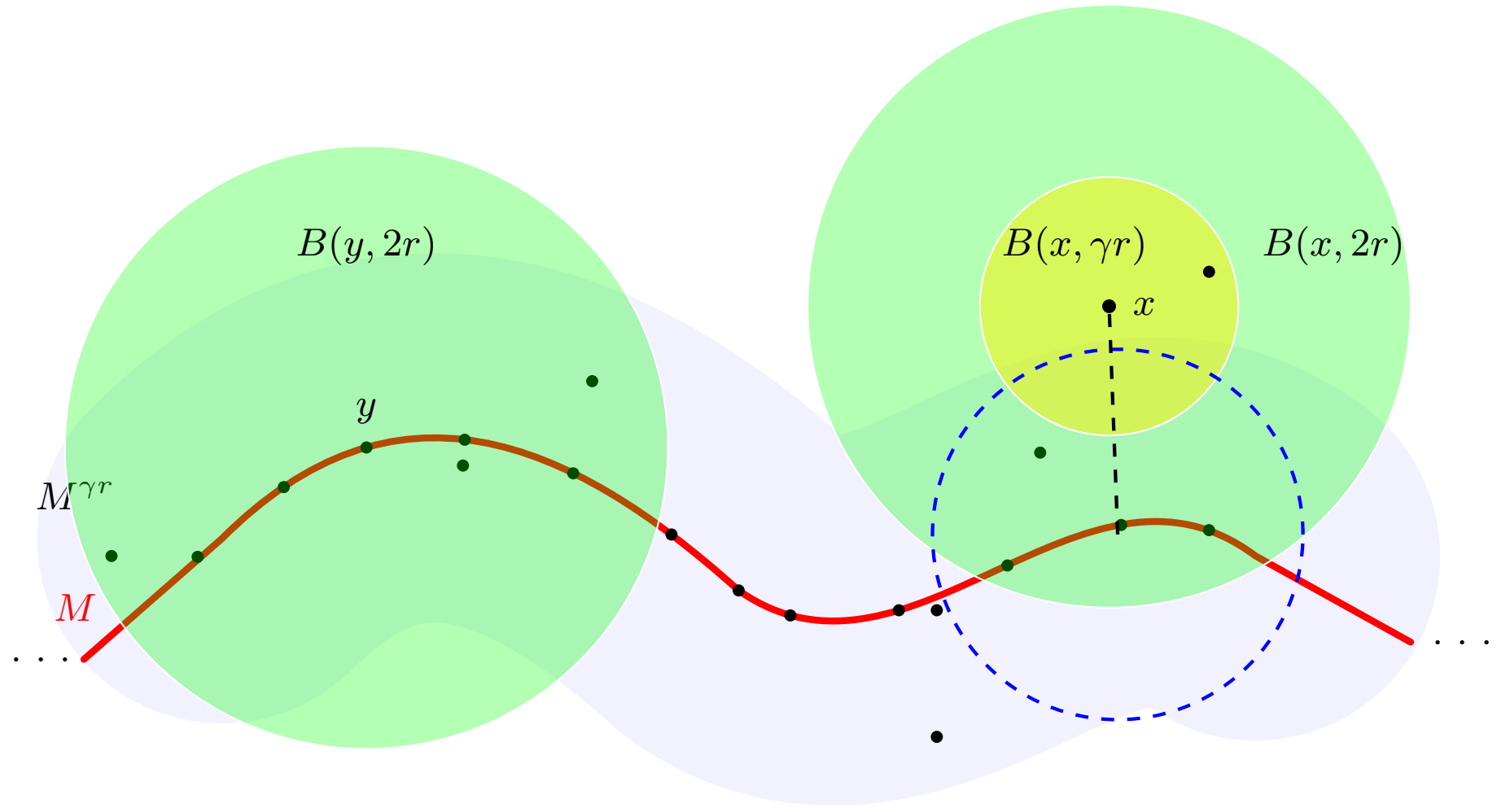

\begin{restatable}{restatablethm}{MainTheoremTwo}\label{thm:hard_prob_inference_bifiltration}
Let $M$ be a $m$-dimensional manifold without boundary with positive reach $\tau$, embedded in the unit cube $[0,1]^d$. 
Let $P$ be a $w$-sample of $M$.
Fix $\delta \in (0,1)$ and suppose that 
$
n \geq \frac{20040^m\Vol(M)}{(1-w)v_m\tau^m }\left( \log\left(\frac{40020^m\Vol(M)}{v_m\tau^m }\right) + \log \frac{2}{\delta}\right)
$.
If 
\begin{align*} r \in\left[\frac{1}{150}\tau, \frac{1}{60}\tau\right]
\end{align*}
and $k \in \mathbb{N}$ satisfies 
\begin{align*}
 k\in \biggl( 1+wnv_d(2r)^d+ \frac{(1-w)nv_m}{\Vol(M)}\left(\frac{91r}{50}\right)^m + \sqrt{\frac{n}{2}\log\frac{2n}{\delta}}, \\ \frac{(1-w)nv_m}{\Vol(M)}\left(\frac{47r}{24}\right)^m-\sqrt{\frac{(1-w)n}{2}\log\frac{2n}{\delta}} \biggr]
\end{align*}
then, the degree-\Cech bifiltration $C_{k,r}(P)$ has the homotopy type of $M$, with probability at least $1-\delta$.
\end{restatable}

The main difference to the first theorem is the additional term in the lower bound for the admissible degree $k$
which comes from the fact that far outliers may have inlier neighbors.
In both theorems, the admissible $k$-interval may be empty for specific parameter choices. We also analyze the special case of a $2$-manifold in $\R^3$ (i.e., a surface)
and show that for a wide range of the signal-to-noise ratio $w$, the $k$-interval will be non-empty for sufficiently large samples.

\subparagraph{Differences to Balakrishnan et al.~\cite{balakrishnan-minimax}}
As remarked in the Related Work Section, our proof follows the
degree-based cleaning strategy of Balakrishnan et al.~\cite{balakrishnan-minimax}:
Their approach in Section 5.2 to get an upper bound is along the same lines as
our approach, filtering out points far away from $M$ and keeping all points
sufficiently close to $M$. Still, their approach does not directly 
yield an inference theorem for the degree \Cech complex:
in their Lemma~3, for instance,
the search radius to determine the neighbors of a point is $2s$, while
the scale of the \Cech complex is $s+\frac{\tau}{2}$ which by the
constraint that $s<(\sqrt{9}-\sqrt{8})\frac{\tau}{2}$ is at least $6.82s$.
Consequently, that result can at most yield inference results
for $c$-degree \Cech complexes with $c<0.15$.
Allowing the full interval of reconstruction radii in Proposition 7.1 of Niyogi et al.~\cite{nsw-finding} extends this range to \(c<1/4\), when using their bounds of \(2s\) on both the distance of retained points to \(M\) and the coverage error.
We obtain guarantees for larger values of \(c\) by employing the reconstruction theorem of Attali et al., which treats the coverage error and the distance of retained points to \(M\) separately (discussed in the next section).
For the standard degree-\Cech complex, corresponding to \(c=1\), we additionally control the number of inlier neighbors of outliers that must be excluded.

\subparagraph{Outline of the rest of the paper}
We first review the result of Attali et al~\cite{akfflsw-tight}
that plays a distinguished role in our approach in Section~\ref{sec:inference}.
We then proceed in Section~\ref{sec:neighbor_count} with counting the
number of neighbors of sample points, subject to the choice $c<\frac{1}{2}$
which leads to the first main theorem.
In Section~\ref{sec:hard_neighbor_count}, we extend the analysis
to the case $c=1$ to
prove the second main theorem.

\section{Topological inference}
\label{sec:inference}
In this section, we establish the geometric ingredients used for the homotopy reconstruction arguments for our results. 
The first is a homotopical inference theorem of Attali et al.~\cite{akfflsw-tight}, which gives sufficient conditions for a union of balls around a point sample to recover the homotopy type of the manifold. 
The second is a coverage condition for points sampled uniformly from the manifold~\cite{nsw-finding}. 


\subparagraph{Homotopical inference theorem of Attali et al.}
Let $X$ be a subset of $\mathbb{R}^d$. The \textit{$r$-thickening} or \textit{$r$-offset} $X^{r}$ of $X$ of radius $r$ is the set of all points of $\mathbb{R}^d$ with a distance at most $r$ from $X$. This is equivalently the union of the ball of radius $r$ centered at the points of $X$, that is, $X^{r} = \bigcup_{x \in X} B(x, r)$. See Figure~\ref{fig:cech_illust}.  

Let $M$ be a compact submanifold of $\mathbb{R}^d$ without boundary.
The \textit{reach} $\tau$ of the manifold $M$ is the supremum of all values $r$ where the nearest projection map  $\pi_0: M^{r} \to M$ defined by $\pi_0(y):=\arg \min_{x \in M} ||x-y||$ is well defined.
Equivalently, the reach is the distance of $M$ to the \emph{medial axis} of $M$, which
is the set of points in the Euclidean space $\mathbb{R}^d$ with more than one nearest points in the manifold $M$.
The reach is a central concept in geometric inference; see \cite{nsw-finding,federer-curvature,berenfeld-estimating} for in-depth studies
of its properties.

We assume that the reader is familiar with basic topological terms such as deformation retract and homotopy equivalence~\cite{h-algebraic}.
Underlying our result is the following homotopy reconstruction of the manifold $M$ from the union of balls $P^r$ of a point set $P$.

\begin{theorem}[Theorem 3, Attali et al.~\cite{akfflsw-tight}] \label{thm:attali_theorem}
    Let $P \subseteq \mathbb{R}^d$ such that $P \subseteq M^{\delta}$ and $M \subseteq P^{\epsilon}$ for $\delta, \epsilon < \tau$.
    Assume that a parameter $\alpha>0$ is small enough, so that the $\alpha$-offset $M^{\alpha}$ of the manifold $M$ is contained in the union of balls $P^r$, i.e. $M^{\alpha} \subseteq P^r$.
    If $r^2 \leq (\tau-\delta)^2 - (\tau-\alpha)^2$, 
    then the union of balls $P^r$ deformation retracts onto $M$. 
\end{theorem}

A major idea of that Theorem is to separate the sampling density $\epsilon$ that controls how well the sample covers the manifold, and the sample noisiness $\delta$ denoting how far
the sample points are allowed to deviate from $M$.
For our work, we will use the following form of inference as a corollary of Theorem~\ref{thm:attali_theorem}.

\begin{restatable}{restatablelemma}{LemmaInfer}\label{lem:infer_thm_r_scale}
     Let $M$ be a submanifold of $\mathbb{R}^d$ with positive reach $\tau$.
     Let $P \subset \mathbb{R}^d$ and $\gamma \in [0,1)$. Set $\epsilon=\frac{(1-\gamma)^2}{4}\tau$ and $\mathcal{R}_{\gamma}=\left[\frac{1-\gamma}{3}\tau, \frac{5(1-\gamma)}{6}\tau\right]$.
     For every $r \in \mathcal{R}_{\gamma}$, if $M \subseteq P^{\epsilon}$ and $P \subseteq M^{\gamma r}$ then, the union of balls $P^r$ deformation retracts onto $M$.
\end{restatable}
\begin{proof}
    Let $\epsilon <\tau$, $\gamma \in (0,1)$ and $r>0$. 
    We take $\alpha=r-\epsilon$, since $M \subseteq P^{\epsilon}$, then we have $M^{\alpha}=M^{r-\epsilon} \subseteq P^{r}$.
    Set $\delta=\gamma r$. The required radius condition of Theorem~\ref{thm:attali_theorem} is  
    \begin{align*}
    r^2 \leq (\tau-\delta)^2-(\tau-\alpha)^2 =  (\tau-\gamma r)^2-(\tau-r+\epsilon)^2.
    \end{align*}
    That is equivalent to
    \begin{align}\label{al:quad_1}
    (2-\gamma^2)r^2 - 2(\tau(1-\gamma)+\epsilon)r+\epsilon(2\tau+\epsilon) \leq 0.
    \end{align}
    One can check that if $\epsilon=\frac{(1-\gamma)^2 \tau}{4}$ and the radius $r$ belongs to the interval $\mathcal{R}_{\gamma}=\left[\frac{1-\gamma}{3}\tau, \frac{5(1-\gamma)}{6}\tau\right]$, then equation (\ref{al:quad_1}) holds. We refer to Appendix \ref{app:interval_recons} for the details.
    Finally, since $P \subseteq M^{\gamma r}$ for $r \in \mathcal{R}_{\gamma}$, Theorem~\ref{thm:attali_theorem} implies that $P^r$ deformation retracts onto $M$.
\end{proof}

We are going to use that lemma in the following form: if a sample $X$ of points in $M$ satisfies $M\subseteq X^\eps$, then we can add
to $X$ an arbitrary number of points within $M^{\gamma r}$ without losing the property that $X^r$ deformation retracts onto $M$.
Importantly for us, the noisiness $\gamma r$ can be much larger than $\eps$; for instance, for $\gamma=\frac{3}{4}$, we get $\epsilon=\frac{\tau}{64}$,
but $\gamma r\geq \frac{3}{4}\cdot\frac{\tau}{12}=\frac{\tau}{16}$. This gap will be crucial in the proof of Theorem~\ref{thm:hard_prob_inference_bifiltration}.

\subparagraph{Manifold coverage.}
In Lemma~\ref{lem:infer_thm_r_scale}, we require that $M\subseteq P^{\epsilon}$. As we will sample $P$
uniformly at random from $M$, we will require a quantitative result how many sample points in $P$
will suffice to ensure that property. This is a direct consequence of a result by Niyogi, Smale, and Weinberger~\cite{nsw-finding}
that we will derive next.
\begin{lemma} \label{lem:epsilon_dense_condition} 
    Let $M$ be a compact submanifold of $\mathbb{R}^d$ of dimension $m$ with positive reach $\tau$.
    Let $\gamma \in (0,1)$ and $\epsilon=\frac{(1-\gamma)^2}{4}\tau$.
    Let $P=\{ x_1,...,x_s\}$ be a set of $s$ points uniformly sampled on $M$ and  $\delta \in (0,1)$. 
    If 
     \[
    s \geq \frac{8.016^m\Vol(M)}{v_m(1-\gamma)^{2m}\tau^m }\left( \log\left(\frac{16.008^m\Vol(M)}{v_m(1-\gamma)^{2m}\tau^m }\right) + \log \frac{2}{\delta}\right)
    \]
    then with probability at least
    $1-\frac{\delta}{2}$, we have $M \subseteq P^{\epsilon}$.
\end{lemma}
\begin{proof}
We employ Prop. 3.2 of Niyogi et al.~\cite{nsw-finding}, which states (up to renaming) that given $\xi<\tau/2$,
given more than
\[\beta_1(\log\beta_2+\log\frac{1}{\delta})\]
sample points with
\[\beta_1=\frac{\Vol(M)}{\cos^m(\arcsin(\frac{\xi}{8\tau}))v_m\left(\frac{\xi}{4}\right)^m},\quad\beta_2=\frac{\Vol(M)}{\cos^m(\arcsin(\frac{\xi}{16\tau}))v_m\left(\frac{\xi}{8}\right)^m}\]
yields a point sample that satifies with probability $1-\delta$ that $M\subseteq P^{\xi/2}$.

We use that proposition in the following way: as we require $M\subseteq P^\epsilon$, we can set $\xi=2\epsilon$ in the bound~-- then, the condition on $\xi$ becomes $\epsilon<\tau/4$ which is satisfied.

For $\beta_1$, setting $\xi=2\epsilon=\frac{(1-\gamma)^2\tau}{2}$ yields
\[\beta_1=\frac{\Vol(M)8^m}{\cos^m( \arcsin(\frac{(1-\gamma)^2}{16}))v_m(1-\gamma)^{2m}\tau^m}.\]

It remains to bound the cosine in the denominator. Note that $\cos(\arcsin(x))=\sqrt{1-x^2}$, hence
\[\cos( \arcsin\left(\frac{(1-\gamma)^2}{16}\right))=\sqrt{1-\left(\frac{(1-\gamma)^2}{16}\right)^2}\geq \sqrt{1-\frac{1}{256}}\]
and since $\frac{8}{\sqrt{1-\frac{1}{256}}}< 8.016$, it follows that $\frac{8.016^m\Vol(M)}{v_m(1-\gamma)^{2m}\tau^m }> \beta_1$.

A completely analogous estimate for $\beta_2$ yields that $\left(\frac{16.008^m\Vol(M)}{v_m(1-\gamma)^{2m}\tau^m }\right)> \beta_2$.
To increase the success probability to $1-\frac{\delta}{2}$, it clearly suffices to replace $\delta$ with $\frac{\delta}{2}$,
and the result follows.
\end{proof}


\section{Proof of~Theorem~\ref{thm:prob_inference_bifiltration}} \label{sec:neighbor_count}
Recall that we refer to $P$ as a $w$-sample of $M$, and that
an \emph{inlier} is a point sampled on $M$, and a \emph{far outlier} is a sample point $x\in P$ with $d(x,M)>2cr$ (with respect to a fixed radius~$r$).

We first fix a far outlier $x\in P$ and count the number of neighbors.
By definition, its neighbors are all within distance $\leq 2cr$,
so no inlier can be a neighbor. On the other hand, another outlier $y$
is a neighbor of $x$ if and only if $y$ lies in the $2cr$-ball around $x$,
and because the outliers are i.i.d. sampled in $[0,1]^d$, the probability
of this event is at most $v_d(2cr)^d$,
where $v_d$ is the volume of the unit $d$-dimensional ball. Because
there are $wn$ outliers drawn, the expected number of neighbors of $x$
is at most
\[\Eout:=1+wnv_d(2cr)^d\]
(where the $1$ comes from counting $x$ as neighbor to itself).
Because this is a Bernoulli sampling, there is a sharp concentration
around the expected value. A simple application of Hoeffding's inequality (Appendix~\ref{app:deg_concentration}) shows that
\[\mathrm{Prob}(\text{$x$ has more than $1+wnv_d(2cr)^d+\sqrt{\frac{wn}{2}\log\frac{2n}{\delta}}$ neighbors})\leq \frac{\delta}{2n}.\]

Next, we fix an inlier $p\in P$ and we employ a similar analysis
except that we want to show a \emph{lower} bound on the number of neighbors.
For that, it suffices to count the expected number of inlier neighbors
of $p$, ignoring potential outliers sampled close to $p$. We use the
following estimate:

\begin{lemma}[Lemma 5.3, Niyogi et al.~\cite{nsw-finding}]\label{lem:lower_intersect_ball_man}
    Let $M$ be a $m$-dimensional submanifold of $\mathbb{R}^d$  with positive reach $\tau$. 
    If $x \in M$, then 
    \[
    \Vol(B(x,s) \cap M) \geq \left(1-\left(\frac{s}{2\tau}\right)^2\right)^{\frac{m}{2}}v_m s^m.
    \]
\end{lemma}

Observe that $q$ is a neighbor of $p$ implies that $q\in B(p,2cr)$,
and because the inliers are sampled i.i.d. from $M$, the probability of this
to happen is at least $\frac{\Vol(B(p,2cr) \cap M)}{\Vol M}$. Using the previous
lemma with $s=2cr$, we obtain that the probability of $q$ being a neighbor
of $p$ is at least
\[
\left(1-\left(\frac{cr}{\tau}\right)^2\right)^{\frac{m}{2}} \frac{v_m (2cr)^m}{\Vol(M)}. 
\]
We impose the condition that $r\leq\frac{5(1-2c)}{6}\tau$ (compare Lemma~\ref{lem:infer_thm_r_scale} with $\gamma=2c$).
Furthermore, for $c\in (0,\frac{1}{2}]$, we have $c(1-2c)\leq \frac{1}{8}$. Therefore,
\[\left(1-\left(\frac{cr}{\tau}\right)^2\right)^{\frac{m}{2}}\geq \left(\sqrt{1-\left(\frac{5}{48}\right)^2}\right)^m \geq \left(\frac{47}{48}\right)^m\]
and the probability that $q$ is neighbor of $p$ is at least
\[\frac{v_m}{\Vol(M)}\left(\frac{47cr}{24}\right)^m.\]
Since the number of inliers is $(1-w)n$, the expected number of neighbors
of $p$ is at least
\[\Ein:=(1-w)n\frac{v_m}{\Vol(M)}\left(\frac{47cr}{24}\right)^m\]
and because the sampling is Bernoulli, Hoeffding's inequality once again yields that
\[\mathrm{Prob}(\text{$p$ has less than $(1-w)n\frac{v_m}{\Vol(M)}\left(\frac{47cr}{24}\right)^m-\sqrt{\frac{(1-w)n}{2}\log\frac{2n}{\delta}}$ neighbors})\leq \frac{\delta}{2n}.\]

Now, a simple invocation of the union bound (See Appendix~\ref{app:deg_concentration}) yields the following intermediate result.

\begin{lemma}\label{lem:prob_all_points_behave_as_they_should}
  Given a $w$-sample $P$ with $n$ points, $\delta\in(0,1)$, $c\in(0,\frac{1}{2})$
  and $r\leq\frac{5(1-2c)}{6}\tau$. With probability at least $1-\frac{\delta}{2}$,
  every far outlier point of $P$ has at most $1+wnv_d(2cr)^d+\sqrt{\frac{wn}{2}\log\frac{2n}{\delta}}$ neighbors and
  every inlier point of $P$ has at least $(1-w)n\frac{v_m}{\Vol(M)}\left(\frac{47cr}{24}\right)^m-\sqrt{\frac{(1-w)n}{2}\log\frac{2n}{\delta}}$ neighbors
  in distance $2cr$.
\end{lemma}

With that lemma, we can prove our first main theorem.
\MainTheoremOne*

\begin{proof}
  Recall that the $c$-degree \Cech complex at $(k,r)$ is just the \Cech complex
  of the point set $P_{k,cr}$ at radius $r$ which is homotopically equivalent to $P_{k,cr}^r$ by the Nerve Theorem. It therefore
  suffices to show that $P_{k,cr}^r$ deformation retracts to $M$.
  
  By Lemma~\ref{lem:prob_all_points_behave_as_they_should}, with probability $1-\frac{\delta}{2}$, no far outlier is in $P_{k,cr}$,
  and every inlier is in $P_{k,cr}$.
  Let $I$ denote the inlier points. By Lemma~\ref{lem:epsilon_dense_condition}
  (with $\gamma=2c$), the inliers satisfy $M\subseteq I^\eps$ with probability $1-\frac{\delta}{2}$, so with probability $1-\delta$, both $I\subset P_{k,cr}$ and
  $M\subseteq I^\eps$.

  Since we ruled out far outliers to be in $P_{k,cr}$, 
  the set $P_{k,cr}\setminus I$ consists of close outliers, all in distance
  at most $2cr$ from $M$. So by definition, $P_{k,cr}\subseteq M^{2cr}$,
  and $M\subseteq P_{k,cr}^{\eps}$ still holds because $I\subseteq P_{k,cr}$.
  Hence with $\gamma=2c$, the conditions of Lemma~\ref{lem:infer_thm_r_scale}
  are satisfied, and $P_{k,cr}^r$ deformation
  retracts to $M$.
\end{proof}

The condition on $k$ in the previous theorem deserves a discussion:
indeed, it is not clear that the given interval for $k$ is non-empty.
Recall $\Ein$ is a lower bound and $\Eout$ is an upper bound
for the expected number of neighbors of an inlier and a far outlier,
respectively, and the bound on $k$ can be rewritten as
\[k\in  \left[\Eout+\sqrt{\frac{wn}{2}\log\frac{2n}{\delta}},\Ein-\sqrt{\frac{(1-w)n}{2}\log\frac{2n}{\delta}} \right]. \]

Obviously, $\Eout<\Ein$ is necessary for the interval to be non-empty.
It is also sufficient in the sense that $\Ein$ and $\Eout$ are both
linear expressions with respect to $n$, so their gap grows linearly,
and the square-root expressions grow asymptotically smaller, so that
for $n$ sufficiently large, the interval will be non-empty
(and also contain an integer value). The precise bound on $n$ depends
on the difference $\Ein-\Eout$, and we omit a detailed analysis.

Instead, we investigate when $\Eout<\Ein$ is satisfied for
the special case of surfaces in three dimensions:

\begin{theorem}\label{thm:special_case_non_empty_interval}
  Let $M$ be a $2$-manifold without boundary in the unit cube $[0,1]^3$.
  For sufficiently many points, $\Eout<\Ein$ holds eventually
  whenever $w<0.483$.
\end{theorem}
\begin{proof}
  Set $A:=v_3(2cr)^3$ and $B:=\frac{v_2}{\Vol(M)}(\frac{47cr}{24})^2$,
  so that $\Eout=1+wnA$ and $\Ein=(1-w)nB$.
  The condition $\Eout<\Ein$ is then equivalent to
  \[
  nB(1-w-w\frac{A}{B})>1.\]
  So if $1-w-w\frac{A}{B}$ is positive, $n$ can be chosen large enough
  so that the inequality is satisfied. This, in turn, is satisfied
  if and only if
  \[w<\frac{1}{1+\frac{A}{B}}\]

  Given that $\frac{v_3}{v_2}=\frac{4}{3}$, an elementary calculation yields
  $\frac{A}{B}\leq 2.782\Vol(M)cr$
  and our upper bound of $cr\leq\frac{5}{6}c(1-2c)\tau\leq \frac{5}{48}\tau$
  gives the bound
  \[\frac{A}{B}\leq 0.29 \Vol(M)\tau.\]
  It remains to bound the quantity $\Vol(M)\tau$. As a consequence
  of Weyl's Tube formula, for any $t<\tau$, the volume of the $t$-thickened
  manifold $M$ is at least (see Lemma~\ref{lem:Weyl} in Appendix~\ref{app:geometric_proofs})
  \[2t(1-\frac{t^2}{3\tau^2})\Vol(M).\]
  On the other hand, that volume is contained in the cube $[-t,1+t]^3$
  of volume $(1+2t)^3$. Choosing $t:=\tau/2$ and exploiting that $\tau<1/2$,
  we obtain that
  \[\Vol(M)\tau\leq\frac{12}{11}(1+\tau)^3.\]
  Finally, we observe that $\tau\leq \frac{1}{2}$: Indeed, consider a maximal
  volume ball inscribed in the interior of the manifold. Such a ball has to touch the manifold in at least two points, as otherwise, the ball could be enlarged.
  Hence, the center of the ball is on the medial axis, and since the ball
  is contained in $[0,1]^3$, its radius is at most $\frac{1}{2}$.
  Hence, $\Vol(M)\tau\leq\frac{81}{22}$
  and therefore, $\frac{A}{B}\leq 1.07$. It follows
  that $\frac{1}{1+\frac{A}{B}}\geq 0.483$ and so,
  if $w$ is smaller than this value, $\Eout<\Ein$ is satisfied for $n$ large enough.
\end{proof}
We point out that the above bound holds uniformly for all values of $c$.
For a fixed $c$, better bounds are possible, especially
if $c$ is close to $\frac{1}{2}$. For instance,
setting $c=0.45$ yields that $cr\leq\frac{5}{6}c(1-2c)\tau=0.0375\tau$
and the same subsequent calculations yield
a non-empty $k$-interval for $w\leq 0.722$.

Similar results can be derived for the case of curves in $[0,1]^3$
(where $\Eout<\Ein$ holds for $w\leq 0.623$) and for curves in $[0,1]^2$
(with $w\leq 0.599$). We give the (slightly simpler) proofs in Appendix~\ref{app:curves}.

\section{Proof of~Theorem~\ref{thm:hard_prob_inference_bifiltration}}\label{sec:hard_neighbor_count}
In the previous section, the condition $c<\frac{1}{2}$ was crucial:
As we set $\gamma:=2c$ in Lemma~\ref{lem:infer_thm_r_scale}, we tolerate
near outliers (with distance $\leq\gamma r$ from $M$) in the sample,
and far outliers (with distance $>\gamma r$) did not have inlier as neighbors,
making counting their (expected) neighbors simple.

For $c\geq\frac{1}{2}$, this analysis fails: Lemma~\ref{lem:infer_thm_r_scale}
requires $\gamma<1$, so we can no longer set $\gamma=2c$. In other words,
no matter what $\gamma$ we chose to define ``near'' and ``far'' outliers,
far outliers can have inliers as neighbors. This leads to a more intricate
analysis of the maximal number of neighbors of far outliers.

While the analysis could be done for arbitrary $c\geq\frac{1}{2}$
and $\gamma\in (0,1)$, we fix henceforth that $c=1$, i.e., we consider
the standard deg-\Cech complex. Moreover, for concreteness,
we set $\gamma:=\frac{49}{50}$. We will comment on the choice at the end of
the paper.
According to the bound of Lemma~\ref{lem:infer_thm_r_scale}, we assume therefore
that $r$ is chosen within the interval $[\frac{\tau}{150},\frac{\tau}{60}]$.
For $r$ fixed in that interval, we call an outlier $y$ \emph{far} if $d(y,M)>\frac{49}{50}r$,
and \emph{near} otherwise.

From now, let $x$ denote a fixed far outlier. We want to upper bound its
expected number of neighbors. There are two types of neighbors: inliers
and outliers. The latter are bounded in the same way as in Section~\ref{sec:neighbor_count}. For the former, we first derive an upper bound on the volume
of the portion of $M$ within distance $2r$ of $x$.
The following lemma states that this intersection is contained in
a ball around the orthogonal projection of $x$ to the manifold.

\begin{restatable}{restatablelemma}{LemmaAnnulus}\label{lem:hard_manifold_annulus}
Let $M \subset \mathbb{R}^d$ be a manifold with positive reach $\tau$.
Let $x \in \mathbb{R}^d$ satisfy $s:=\operatorname{dist}(x,M) < \tau$ and let $z = \pi_0(x)$ be the closest point on $M$ to $x$.
If $s < 2r < \tau$, then
\[
M \cap B(x,2r) \subseteq M \cap B\left(z, \sqrt{\frac{(2r)^2-s^2}{1-\frac{s}{\tau}}}\right).
\]
\end{restatable}

The proof exploits that the manifold is approximately orthogonal to the vector
$x-z$ in an $r$-ball around $x$, provided that $r$ is small enough.
See Appendix~\ref{app:geometric_proofs} for the (elementary) argument.

The second result we need is an upper bound for the volume
of the intersection of $M$ with a (small) ball around a manifold point:

\begin{restatable}{restatablelemma}{LemmaVolUp}\label{lem:hard_vol_upper}
    Let $z \in M$ and let $0 < \rho < \frac{\tau}{2}$. 
    Then 
    \[
    \Vol(M \cap B(z,\rho) ) \leq \frac{v_m \rho^m}{(1-\frac{2 \rho }{\tau})^{\frac{m}{2}}}.
    \]
\end{restatable}

The proof uses a volume distortion argument, measuring how much can the volume
of $M \cap B(z,\rho)$ change when projected to $T_z M\cap B(z,\rho)$,
where $T_z M$ is the tangent space of $M$ at $z$. That requires an investigation
of the Jacobian and several other non-elementary geometric concepts which are
presented in Appendix~\ref{app:geometric_proofs}.

\medskip

We apply these lemmas in our setting next.

\begin{lemma}\label{lem:eout_prime_first_step}
  Let $x$ be a far outlier and $z=\pi_0(x)$. Then
  \[M\cap B(x,2r) \subseteq M\cap B(z,\frac{30}{17}r).\]
\end{lemma}
\begin{proof}
The statement is trivially true if $\|x-z\|>2r$, as the left hand side is the empty set.
Otherwise, we are in the situation of Lemma~\ref{lem:hard_manifold_annulus}.
The term $\frac{(2r)^2-s^2}{1-\frac{s}{\tau}}$ is monotonously decreasing in $s$ in the interval $[\frac{49}{50}r,2r]$.
Hence, the maximum radius is obtained for $s=\frac{49}{50}r$ and we calculate
\begin{align*}
    \frac{4r^2 - (\frac{49}{50}r)^2}{1-\frac{49r}{50\tau}}\leq r^2\frac{45594}{14755} < r^2 \left(\frac{30}{17}\right)^2,
\end{align*}
where the first inequality follows from $r \leq \frac{1}{60}\tau$.
\end{proof}

\begin{lemma}
  Let $x$ be a far outlier. Then, the expected number of inlier neighbors of $x$ is at most
  \[
  \Eout' := (1-w)n\frac{v_m}{\Vol(M)}\left(\frac{91r}{50}\right)^m.
  \]
\end{lemma}
\begin{proof}
  Applying Lemma~\ref{lem:eout_prime_first_step} and Lemma~\ref{lem:hard_vol_upper} with $\rho=\frac{30r}{17}$ we get
\[
\Vol(M \cap B(x,2r)) 
\leq v_m \left(\frac{30r}{17}\right)^m \left(1-\frac{60r}{17\tau}\right)^{-\frac{m}{2}}.
\]
We observe that $\sqrt{1-\frac{60r}{17\tau}} \geq \sqrt{1-\frac{1}{17}}=\sqrt{\frac{16}{17}}\geq \frac{32}{33}$, then the volume above simplifies to
\[
\Vol(M \cap B(x,2r)) \leq v_m \left(\frac{30r}{17}\right)^m\left(\frac{33}{32}\right)^m\leq v_m\left(\frac{91r}{50}\right)^m.
\]
The claim follows because there are $(1-w)n$ inliers, and each of them is neighbor of $x$ with probability at most $\frac{\Vol(M \cap B(x,2r))}{\Vol(M)}$.
\end{proof}

Since a far outlier may receive neighbors both outliers and inliers,
its expected number of neighbors is at most 
\[\Eout + \Eout'\]
with $\Eout$ as in Section~\ref{sec:neighbor_count}.
Again, applying the Hoeffding's inequality (See Appendix~\ref{app:deg_concentration}), we have  
\[
\mathrm{Prob}(\text{$x$ has more than $1+wnv_d(2cr)^d+\frac{(1-w)nv_m}{\Vol(M)}\left(\frac{91r}{50}\right)^m+\sqrt{\frac{n}{2}\log\frac{2n}{\delta}}$ neighbors})\leq \frac{\delta}{2n}.
\]

The lower bound for the expected number of neighbors of an inlier point remains $\Ein$, as in Section~\ref{sec:neighbor_count}.
We are ready to state the main reconstruction theorem for the case of the degree-\Cech bifiltration.

\MainTheoremTwo*
\begin{proof}
  The proof proceeds in a similar fashion as for Theorem~\ref{thm:prob_inference_bifiltration}.
  The bounds for $n$ and $r$ come from Lemma~\ref{lem:infer_thm_r_scale} and Lemma~\ref{lem:epsilon_dense_condition} by the choice $\gamma=\frac{49}{50}$.
  The $k$-interval is the consequence of the expected number $\Ein$ of neighbors of inlier points, and the expected number $\Eout + \Eout'$ of far outlier
  points as analyzed above: if $k$ is in the interval, with probability $1-\frac{\delta}{2}$, all inliers are in $P_{k,r}$, and no far outlier is.
  Moreover, with same probability, the inliers form an $\epsilon$-dense sample of $M$, and Lemma~\ref{lem:infer_thm_r_scale}
  yields homotopy equivalence with the degree-\Cech complex.
\end{proof}

Again, the question remains whether the $k$-interval is non-empty
(for sufficiently many points). Note that $\frac{91}{50}<\frac{47}{24}$,
which is necessary. Similar as for Theorem~\ref{thm:prob_inference_bifiltration},
it is enough to analyze when the expected values satisfy
$\Eout+\Eout'<\Ein$. The analysis for $2$-manifolds in $\R^3$
is along the same lines as for Theorem~\ref{thm:special_case_non_empty_interval}:

\begin{theorem}\label{thm:special_case_non_empty_interval_deg_cech}
  Let $M$ be a $2$-manifold without boundary in the unit cube $[0,1]^3$.
  For sufficiently many points, $\Eout+\Eout'<\Ein$ holds eventually
  whenever $w<0.443$.
\end{theorem}
\begin{proof}
  Set $A':=v_3(2r)^3$ and $B':=\frac{v_2r^2}{\Vol(M)}(\left(\frac{47}{24}\right)^2-\left(\frac{91}{50}\right)^2)$. With that, we have
  \[\Ein-(\Eout+\Eout')=nB'(1-w-w\frac{A'}{B'})-1,\]
  which means that for $n$ large enough, $\Eout+\Eout'<\Ein$ will be
  satisfied if $w<\frac{1}{1+\frac{A'}{B'}}$.
  We can calculate that
  \[\frac{A'}{B'}=\frac{32}{3(\left(\frac{47}{24}\right)^2-\left(\frac{91}{50}\right)^2)}\Vol(M)r\leq 20.41 \Vol(M)r.\]
  Now using $r\leq\frac{\tau}{60}$ and $\Vol(M)\tau\leq\frac{81}{22}$
  as in Theorem~\ref{thm:special_case_non_empty_interval} yields
  $\frac{A'}{B'}\leq 1.253$
  and so
  \[\frac{1}{1+\frac{A'}{B'}}\geq 0.443.\qedhere\]
\end{proof}
  
\section{Discussion and outlook}\label{sec:conclusion}
\subparagraph{Sample sizes.} An obvious drawback of Theorem~\ref{thm:hard_prob_inference_bifiltration}
is a large sample size required to provide a guarantee.
For instance, fixing $M$ to be a sphere in $\R^3$ of radius $0.4$, setting $w=0.1$ and $\delta=0.1$
yields a sample size of at least $45.6$ billion points. The reason for that large size is
of course the choice of $\gamma$ as $\frac{49}{50}$ which causes $\epsilon$ to be $\frac{\tau}{10000}$.
On the other hand, a smaller value of $\gamma$ leads to an increase of the basis $\frac{91r}{50}$ in
Theorem~\ref{thm:hard_prob_inference_bifiltration} which has to stay below $\frac{47r}{24}$.

We believe that a more careful analysis could lead to significantly smaller sample sizes~--
we made little attempt to keep the constants in our bounds small. We also speculate that
a tailored analysis towards low-dimensional point samples would give improved bounds for the
special case.
Finally, we remark that the sample size in Theorem~\ref{thm:prob_inference_bifiltration}
is more reasonable, especially for moderate values of $c$. For instance, for the same sphere as above
and choosing $c=\frac{1}{4}$, reconstruction is guaranteed for a choice of $r$ and $k$ with only
around $150,000$ points.

\subparagraph{Further work.}
We consider our work as a case study to demonstrate that probabilistic inference results
are possible also in the presence of outlier noise, in the spirit
of the result of Balakrishnan et al.~\cite{balakrishnan-minimax}.
We believe that such results
strengthen the argument
that bifiltrations such as the degree-\Cech are promising tools for data analysis in the
presence of such outliers. There are many scenarios which can be studied next.
For instance, one could consider degree-Rips complexes instead, change to a
different noise model (e.g., where outliers are sampled normally according to the distance
to the manifold) or consider the case of manifolds with boundary.

\bibliographystyle{plainurl}
\bibliography{refs}

@book{bcy-geometric,
  title     = {Geometric and Topological Inference},
  author    = {Boissonnat, Jean-Daniel and Chazal, Fr{\'e}d{\'e}ric and Yvinec, Mariette},
  series    = {Cambridge Texts in Applied Mathematics},
  year      = {2018},
  publisher = {Cambridge University Press},
  doi       = {10.1017/9781108226653}
}

@article{mz-manifold,
   author = "Meilă, Marina and Zhang, Hanyu",
   title = "Manifold Learning: What, How, and Why", 
   journal= "Annual Review of Statistics and Its Application",
   year = 2024,
   volume = 11,
   pages = "393-417",
   doi = "https://doi.org/10.1146/annurev-statistics-040522-115238",
  }

@article{nsw-finding,
  author    = {Niyogi, Partha and Smale, Stephen and Weinberger, Shmuel},
  title     = {Finding the Homology of Submanifolds with High Confidence from Random Samples},
  journal   = {Discrete \& Computational Geometry},
  volume    = 39,
  number    = 1,
  pages     = {419--441},
  year      = 2008,
  publisher = {Springer},
  doi       = {10.1007/s00454-008-9053-2}
}

@article{federer-curvature,
 ISSN = {00029947},
 URL = {http://www.jstor.org/stable/1993504},
 author = {Herbert Federer},
 journal = {Transactions of the American Mathematical Society},
 number = {3},
 pages = {418--491},
 publisher = {American Mathematical Society},
 title = {Curvature Measures},
 urldate = {2026-06-16},
 volume = {93},
 year = {1959}
}

@article{berenfeld-estimating,
  title={Estimating the reach of a manifold via its convexity defect function},
  author={Berenfeld, Cl{\'e}ment and Harvey, John and Hoffmann, Marc and Shankar, Krishnan},
  journal={Discrete \& Computational Geometry},
  volume={67},
  number={2},
  pages={403--438},
  year={2022},
  publisher={Springer}
}

@article{bl-l-stability,
  author       = {Andrew J. Blumberg and
                  Michael Lesnick},
  title        = {Stability of 2-Parameter Persistent Homology},
  journal      = {Found. Comput. Math.},
  volume       = {24},
  number       = {2},
  pages        = {385--427},
  year         = {2024},
  url          = {https://doi.org/10.1007/s10208-022-09576-6},
  doi          = {10.1007/S10208-022-09576-6},
}

@article{r-s-stable-clustering,
  author  = {Alexander Rolle and Luis Scoccola},
  title   = {Stable and Consistent Density-Based Clustering via Multiparameter Persistence},
  journal = {Journal of Machine Learning Research},
  year    = {2024},
  volume  = {25},
  number  = {258},
  pages   = {1--74},
  url     = {http://jmlr.org/papers/v25/21-1185.html}
}

@incollection{bl-introduction,
  author    = {Botnan, Magnus Bakke and Lesnick, Michael},
  title     = {An Introduction to Multiparameter Persistence},
  booktitle = {Proceedings of the 2020 International Conference on Representations of Algebras},
  pages     = {77--150},
  year      = 2023,
  publisher = {EMS Press},
  doi       = {10.4171/282}
}

@inproceedings{akfflsw-tight,
  author    = {Attali, Dominique and Dal Poz Kou{\v{r}}imsk{\'a}, Hana and Fillmore, Christopher and Ghosh, Ishika and Lieutier, Andr{\'e} and Stephenson, Elizabeth and Wintraecken, Mathijs},
  title     = {Tight Bounds for the Learning of Homotopy {\`{a}} la Niyogi, Smale, and Weinberger for Subsets of Euclidean Spaces and of Riemannian Manifolds},
  booktitle = {40th International Symposium on Computational Geometry (SoCG 2024)},
  volume    = 293,
  pages     = {11:1--11:19},
  year      = 2024,
  doi       = {10.4230/LIPIcs.SoCG.2024.11},
}

@article{bbgs-core,
  author  = {Nello Blaser and Morten Brun and Odin Hoff Gard{\aa} and Lars M. Salbu},
  title   = {Core bifiltration},
  journal = {Journal of Applied and Computational Topology},
  volume  = 9,
  number  = 4,
  pages   = 30,
  year    = 2025,
  month   = {November},
  doi     = {10.1007/s41468-025-00226-8},
}

@book{h-algebraic,
  title={Algebraic Topology},
  author={Hatcher, Allen},
  year={2001},
  publisher={Cambridge University Press},
  address={Cambridge},
  isbn={9780521795401}
}

@InProceedings{a-role-deg-rips-annulus,
  author =  {Rolle, Alexander},
  title = {{The Degree-Rips Complexes of an Annulus with Outliers}},
  booktitle = {38th International Symposium on Computational Geometry (SoCG 2022)},
  pages = {58:1--58:14},
  series =  {Leibniz International Proceedings in Informatics (LIPIcs)},
  ISBN =  {978-3-95977-227-3},
  year =  {2022},
  volume =  {224},
  doi =   {10.4230/LIPIcs.SoCG.2022.58},
}

@article{c-l-manifol-reconst,
title = {Smooth manifold reconstruction from noisy and non-uniform approximation with guarantees},
journal = {Computational Geometry},
volume = {40},
number = {2},
pages = {156-170},
year = {2008},
issn = {0925-7721},
doi = {https://doi.org/10.1016/j.comgeo.2007.07.001},
author = {Frédéric Chazal and André Lieutier},
}

@incollection{h-vietoris-rips,
  author    = {Hausmann, Jean-Claude},
  title     = {On the {V}ietoris-{R}ips complexes and a cohomology theory for metric spaces},
  booktitle = {Prospects in Topology: Proceedings of a Conference in Honor of William Browder},
  editor    = {Quinn, Frank},
  series    = {Annals of Mathematics Studies},
  volume    = {138},
  pages     = {177--188},
  year      = {1994},
  publisher = {Princeton University Press},
  address   = {Princeton, NJ}
}

@inproceedings{c-o-towards-pers-recons,
  author    = {Fr\a'ed\a'eric Chazal and Steve Yann Oudot},
  title     = {Towards persistence-based reconstruction in Euclidean spaces},
  booktitle = {Proceedings of the twenty-fourth annual symposium on Computational geometry},
  pages     = {232--241},
  year      = {2008},
  publisher = {ACM},
  doi       = {10.1145/1377676.1377719}
}

@InProceedings{k-j-c-w-homot-recons,
  author =  {Kim, Jisu and Shin, Jaehyeok and Chazal, Fr\'{e}d\'{e}ric and Rinaldo, Alessandro and Wasserman, Larry},
  title = {{Homotopy Reconstruction via the Cech Complex and the Vietoris-Rips Complex}},
  booktitle = {36th International Symposium on Computational Geometry (SoCG 2020)},
  pages = {54:1--54:19},
  series =  {Leibniz International Proceedings in Informatics (LIPIcs)},
  ISBN =  {978-3-95977-143-6},
  ISSN =  {1868-8969},
  year =  {2020},
  volume =  {164},
  doi =   {10.4230/LIPIcs.SoCG.2020.54},

}

@article{c-s-q-geom-infer-prob-meas,
  author    = {Chazal, Fr{\'e}d{\'e}ric and Cohen-Steiner, David and M{\'e}rigot, Quentin},
  title     = {Geometric Inference for Probability Measures},
  journal   = {Foundations of Computational Mathematics},
  volume    = {11},
  number    = {4},
  pages     = {433--451},
  year      = {2011},
  publisher = {Springer},
  doi       = {10.1007/s10208-011-9098-0},
  url       = {https://link.springer.com/article/10.1007/s10208-011-9098-0}
}

@article{n-s-w-top-learn,
  author    = {Niyogi, Partha and Smale, Steve and Weinberger, Shmuel},
  title     = {A Topological View of Unsupervised Learning from Noisy Data},
  journal   = {SIAM Journal on Computing},
  volume    = {40},
  number    = {3},
  pages     = {646--663},
  year      = {2011},
  publisher = {Society for Industrial and Applied Mathematics},
  doi       = {10.1137/090762932},
}

@article{blw-reach,
  author  = {Boissonnat, Jean-Daniel and Lieutier, Andr{\'e}
             and Wintraecken, Mathijs},
  title   = {The reach, metric distortion, geodesic convexity
             and the variation of tangent spaces},
  journal = {Journal of Applied and Computational Topology},
  volume  = 3,
  pages   = {29--58},
  year    = 2019,
  doi     = {10.1007/s41468-019-00029-8}
}

@article{acc-minimax,
  author  = {Arias-Castro, Ery and Chau, Phong Alain},
  title   = {Minimax estimation of distances on a surface and minimax
             manifold learning in the isometric-to-convex setting},
  journal = {Information and Inference: A Journal of the IMA},
  volume  = 12,
  number  = 4,
  pages   = {2970--3009},
  year    = 2023,
  doi     = {10.1093/imaiai/iaad046}
}

@misc{Simon2014,
  author = {Simon, Leon},
  title  = {Introduction to Geometric Measure Theory},
  year   = 2014,
  note   = {Lecture notes, Stanford University},
  url    = {https://web.stanford.edu/class/math285/ts-gmt}
}

@misc{wang2019weyl,
  author       = {Wang, Yi},
  title        = {{Weyl}'s Tube Formula, the {Heintze--Karcher}
                  Inequality, and their Consequences in Comparison Geometry},
  year         = 2019,
  howpublished = {Undergraduate thesis},
  note         = {2018--2019; supervised by Simon Brendle},
  url          = {https://bpb-us-w2.wpmucdn.com/web.sas.upenn.edu/dist/6/641/files/2019/09/YiWang_SeniorThesis.pdf}
}

@article{burtscherheckman2021weyl,
  author  = {Burtscher, Annegret and Heckman, Gert},
  title   = {Variations of {Weyl}'s Tube Formula},
  journal = {The Journal of Geometric Analysis},
  year    = 2021,
  volume  = 31,
  pages   = {11952--11970},
  doi     = {10.1007/s12220-021-00707-z},
  url     = {https://doi.org/10.1007/s12220-021-00707-z}
}

@book{frieze_karonski_2016, 
  title={Introduction to Random Graphs}, 
  publisher={Cambridge University Press}, 
  author={Frieze, Alan and Karo{\'n}ski, Micha{\l}}, 
  year={2016},
  doi={10.1017/CBO9781316339831}
}

@article{akcmrw-estimating,
  author  = {Aamari, Eddie and Kim, Jisu and Chazal, Fr{\'e}d{\'e}ric
             and Michel, Bertrand and Rinaldo, Alessandro and Wasserman, Larry},
  title   = {Estimating the Reach of a Manifold},
  journal = {Electronic Journal of Statistics},
  volume  = {13},
  number  = {1},
  pages   = {1359--1399},
  year    = {2019},
  doi     = {10.1214/19-EJS1551}
}

@article{bkrr-unified-nerve,
title = {A unified view on the functorial nerve theorem and its variations},
journal = {Expositiones Mathematicae},
volume = {41},
number = {4},
pages = {125503},
year = {2023},
doi = {https://doi.org/10.1016/j.exmath.2023.04.005},
author = {Ulrich Bauer and Michael Kerber and Fabian Roll and Alexander Rolle},
}

@InProceedings{balakrishnan-minimax,
  title =    {Minimax rates for homology inference},
  author =   {Balakrishnan, Sivaraman and Rinaldo, Alesandro and Sheehy, Don and Singh, Aarti and Wasserman, Larry},
  booktitle =    {Proceedings of the Fifteenth International Conference on Artificial Intelligence and Statistics},
  pages =    {64--72},
  year =   {2012},
  volume =   {22},
  series =   {Proceedings of Machine Learning Research},
}

\appendix
\section{Interval reconstruction}\label{app:interval_recons}
\LemmaInfer*
\begin{proof}[Completion of the proof of Lemma~\ref{lem:infer_thm_r_scale}]
Let $\epsilon <\tau$, $\gamma \in [0,1)$ and $\alpha=r-\epsilon$. 

On the one hand, we show that $M^{\alpha} \subseteq P^{r}$.
If $y \in M^{r-\epsilon}$, then there exists $z \in M$ such that $||y-z||\leq r-\epsilon$. Since $M \subseteq P^{\epsilon}$, there exists $p\in P$ with $||z-p||\leq \epsilon$, and therefore $||y-p||\leq||y-z||+||z-p|| \leq r$. 
Thus $y\in P^r$.

On the other hand, by setting $\delta=\gamma r$ and $\alpha=r-\epsilon$, the required radius condition from Theorem~\ref{thm:attali_theorem} is
\begin{align*}
r^2 \leq (\tau-\delta)^2-(\tau-\alpha)^2= (\tau-\gamma r)^2-(\tau-r+\epsilon)^2.
\end{align*}
Expanding the squares on the right hand side, gives  a quadratic inequality in terms of $r$
\begin{align*}
Q(r) := (2-\gamma^2)r^2 - 2\big(\tau(1-\gamma)+\epsilon\big)r+\epsilon(2\tau+\epsilon) \leq 0.
\end{align*}
Because $\gamma \in [0,1)$, its leading coefficient $(2-\gamma^2)$ is strictly positive, meaning $Q(r)$ defines a convex parabola. 
To prove that $Q(r) \leq 0$ for all $r$ in the interval $\mathcal{R}_{\gamma}=\left[\frac{1-\gamma}{3}\tau, \frac{5(1-\gamma)}{6}\tau\right]$, it is sufficient to show that $Q(r)$ is less than or equal to zero at the two endpoints of the interval.

First, we substitute $\epsilon = \frac{(1-\gamma)^2}{4}\tau$ into $Q(r)$ and evaluate the lower endpoint $r_{\min} = \frac{1-\gamma}{3}\tau$, that yields
\begin{align*}
Q(r_{\min}) &= -\frac{\tau^2(1-\gamma)^2}{144} \Big( 7\gamma^2 - 6\gamma + 7 \Big) 
            = -\frac{\tau^2(1-\gamma)^2}{144} \left[ 7\left(\gamma - \frac{3}{7}\right)^2 + \frac{40}{7} \right].
\end{align*}
The latter is strictly negative, $Q(r_{\min}) < 0$.

Next, we evaluate the upper endpoint $r_{\max} = \frac{5(1-\gamma)}{6}\tau$,

\begin{align*}
Q(r_{\max}) &= -\frac{\tau^2(1-\gamma)^2}{144} \Big( 91\gamma^2 - 42\gamma + 19 \Big)
            = -\frac{\tau^2(1-\gamma)^2}{144} \left[ 91\left(\gamma - \frac{3}{13}\right)^2 + \frac{184}{13} \right].
\end{align*}
Again, this is strictly negative $Q(r_{\max}) < 0$.

Since $Q(r)$ is negative at both endpoints of $\mathcal{R}_{\gamma}$ and is a parabola, $Q(r) \leq 0$ is satisfied for all $r \in \mathcal{R}_{\gamma}$. 
\end{proof}

\section{Degree concentration}\label{app:deg_concentration}
We use the following standard Hoeffding's inequality~\cite[Corollary 30.10]{frieze_karonski_2016}.
Let $X$ be a sum of $N$ independent Bernoulli random variables with mean $\mu$. For every $s>0$, we have
\[
    \mathrm{Prob}\left(X-\mu\geq s\right)
    \leq e^{-\frac{2s^2}{N}},
\text{ and }
    \mathrm{Prob}\left(X-\mu\leq -s\right)
    \leq e^{-\frac{2s^2}{N}}.
\]
\subparagraph{Far outliers has no inlier neighbors.}
Let $\delta\in(0,1)$ and $r>0$.
Recall that
\[
\Eout:=1+wnv_d(2cr)^d
\text{ and }
\Ein:=(1-w)n\frac{v_m}{\Vol(M)}
\left(\frac{47cr}{24}\right)^m.
\]
Let $x\in P$ be a far outlier. We condition on the location of $x$.
Since $x$ is counted as its own neighbor, its neighbor count is $1$ plus a sum of $wn-1$ independent Bernoulli random variables. 
Its conditional expectation is then bounded above by
\[
1+(wn-1)v_d(2cr)^d\leq\Eout.
\]
Let $A_x$ be the event that $x$ is a far outlier and has more than
\[
1+wnv_d(2cr)^d+ \sqrt{\frac{wn}{2}\log\frac{2n}{\delta}}
\]
 neighbors.
Conditional on the location of $x$, we apply Hoeffding's inequality and $wn-1\leq wn$. Then
\[
\mathrm{Prob}(A_x)\leq
\exp\left(-\frac{2}{wn}\left(\sqrt{\frac{wn}{2}\log\frac{2n}{\delta}}\right)^2
\right)
=\frac{\delta}{2n}.
\]

Similarly, let $p\in P\cap M$ be an inlier.
Conditional on the location of $x$. 
We ignore the outlier neighbors.
The neighbor count is bounded below by $1$ plus a sum of $(1-w)n-1$ independent Bernoulli random variables.
Then, the conditional expectation of this lower bound is
\[
1+((1-w)n-1)\frac{\Vol(B(p,2cr)\cap M)}{\Vol(M)} 
\geq (1-w)n\frac{v_m}{\Vol(M)} \left(\frac{47cr}{24}\right)^m
=\Ein.
\]
Let $A_p$ be the event that $p$ has less than
\[
(1-w)n\frac{v_m}{\Vol(M)}
\left(\frac{47cr}{24}\right)^m -\sqrt{\frac{(1-w)n}{2}\log\frac{2n}{\delta}}
\]
neighbors in distance $2cr$.
Again, conditional on $p$, Hoeffding's inequality and $(1-w)n-1\leq(1-w)n$ gives
\[
\mathrm{Prob}(A_p)\leq
\exp\left(-\frac{2}{(1-w)n}\left(\sqrt{\frac{(1-w)n}{2}\log\frac{2n}{\delta}} \right)^2
\right)
=\frac{\delta}{2n}.
\]
Because these bounds are uniform over all locations of the inliers and far outliers, they hold unconditionally.

We now establish the probability that the desired bounds hold simultaneously for sample points in $P$. 
We apply the union bound over all $n$ points in $P$:
\begin{align*}
    \mathrm{Prob}\left( \bigcup_{x} A_x \cup \bigcup_{p } A_p \right) &\leq \sum_{x} \mathrm{Prob}(A_x) + \sum_{p } \mathrm{Prob}(A_p)
    = n \left(\frac{\delta}{2n}\right) = \frac{\delta}{2}.
\end{align*}

Therefore, the probability that none of these events occur (i.e., all far outliers satisfy their upper bound and all inlier satisfy their lower bound )is at least $1-\frac{\delta}{2}$.

\subparagraph{Far outliers have both inlier and outliers neighbors.}
Consider a far outlier $x\in P$ and condition on its location. 
Its neighbors within distance $2r$ are drawn either from the ambient noise or from the manifold,
i.e., from two different Bernoulli distributions.
Excluding $x$ itself, these contributions consist of $wn-1$ and $(1-w)n$ independent Bernoulli random variables, respectively.

Therefore, the conditional expected number of neighbors of $x$ is bounded by $\Eout+\Eout'$, where
\[
\Eout'=(1-w)n\frac{v_m}{\Vol(M)}
\left(\frac{91r}{50}\right)^m.
\]
Although these Bernoulli variables are not identically distributed, they are independent and bounded.
Thus, Hoeffding's inequality applies to their sum.

Let $A_x$ be the event that $x$ is a far outlier and has more than
\[
1+wnv_d(2r)^d+ \frac{(1-w)nv_m}{\Vol(M)} \left(\frac{91r}{50}\right)^m +\sqrt{\frac{n}{2}\log\frac{2n}{\delta}}
\]
neighbors.
Conditional on $x$, there are $n-1$ contributions.
Then, by  Hoeffding's inequality we have
\[
\mathrm{Prob}(A_x)\leq \exp\left(-\frac{2}{n}\left( \sqrt{\frac{n}{2}\log\frac{2n}{\delta}} \right)^2 \right)
=\frac{\delta}{2n}.
\]
Again, this bound is uniform over all far outlier locations and hence holds unconditionally.

\section{Proof of geometric lemmas}\label{app:geometric_proofs}

\begin{lemma}\label{lem:Weyl}
  For $t<\tau$, let $M^t$ denote the $t$-thickening of a $2$-dimensional
  manifold $M$, that is,
  \[M^t:=\{x\in\R^3\mid \operatorname{dist}(x,M)\leq t\}.\]
  Then, we have
  \[\Vol(M^t)\geq 2t(1-\frac{t^2}{3\tau^2})\Vol(M).\]
\end{lemma}
\begin{proof}
  By Weyl's tube formula~\cite[Eq. (52)]{wang2019weyl}, \cite{burtscherheckman2021weyl},
  \[
  \Vol(M^t)=2t\Vol(M)+\frac{2t^3}{3}\int_M \kappa_1(p)\kappa_2(p)dA(p).
  \]
  Moreover, the principal curvatures are lower-bounded by the reach via
  $\kappa_1(p)\kappa_2(p)\geq-\frac{1}{\tau^2}$ (\cite[Lemma 4]{blw-reach}).
  It follows that
  \[\int_M \kappa_1(p)\kappa_2(p)dA(p)\geq -\frac{\Vol(M)}{\tau^2}\]
  and the result follows by plugging into Weyl's formula.
\end{proof}

\LemmaAnnulus*
\begin{proof}
  We can assume that $s>0$ as otherwise $x=z$ and the statement is clear.
  Since $z$ is the unique nearest point to $x$, we may write $x = z + s\nu $, where $\nu$ is a unit normal vector to $M$ at $z$.
    We first claim that for $y\in M$, 
    \begin{equation}
      \label{interm_claim}
    \langle y-z, \nu\rangle  \leq \frac{||y-z||^2}{2\tau}, \ \  y \in M.
    \end{equation}
    Indeed, for $0 < t< \tau$, the point $z$ is the unique nearest point of $M$ to $z + t\nu $~\cite[Thm 4.8 (12)]{federer-curvature}, \cite[Thm 2]{blw-reach}.
    Thus
    \[
    ||y-(z+t \nu)||^2 \geq t^2.
    \]
    Since
    \[||y-(z+t \nu)||^2=||(y-z)-t \nu)||^2=||(y-z)||^2+t^2-2t\langle y-z,\nu\rangle,\]
    it follows $2t \langle y-z, \nu \rangle  \leq ||y -z||^2$ for all $0<t<\tau$. Hence, the inequality also holds in the limit $t= \tau$, proving (\ref{interm_claim}).
    
    Now consider $y\in M\cap B(x,2r)$. Using again that $x=z+s\nu$, we obtain with (\ref{interm_claim}) that
    \begin{align*}
    ||x-y||^2 =||s \nu -(y-z)||^2= &s^2+||y-z||^2-2s\langle \nu , y-z \rangle\\
                                 \geq  & s^2 + \left( 1- \frac{s}{\tau}\right)||y -z||^2.
    \end{align*}
    Since $||x -y|| \leq 2r$ and $s < \tau$, it follows that  
    \[
    ||y-z||^2 \leq \frac{(2r)^2-s^2}{1-\frac{s}{\tau}},
    \]
    and the claim follows by taking square roots.
\end{proof}

\LemmaVolUp*
\begin{proof}
  Let $T_zM$ be the tangent space of $M$ at $z$. Consider the orthogonal projection $\pi: M\to T_zM$. On $M\cap B(z,\rho)$, $\pi$ is injective (Niyogi, Smale, Weinberger
  prove the non-singularity of the projection~\cite[Lemma 5.4]{nsw-finding}; see~\cite[Lemma 2.11]{acc-minimax}
  for a complete proof of injectivity). Write $D:=\pi(M\cap B(z,\rho))$ for the projection. Since
  orthogonal projection does not increase distances, we have that $D$ is contained in $T_zM \cap B(z,\rho)$ which is an $m$-dimensional
  Euclidean ball with volume $v_m \rho^m$.

  It remains to bound the maximal volume distortion of the projection. Fix $y\in M\cap B(z,\rho)$.
  First, the angle $\theta$ 
  between the tangent space $T_zM$ and $T_yM$ satisfies
  \[
    \cos \theta  \geq \sqrt{1-\frac{2||y-z||}{\tau}}\geq \sqrt{1 -\frac{2 \rho}{\tau}}=:\alpha
  \]
  by Niyogi, Smale, Weinberger~\cite[Prop 6.2+6.3]{nsw-finding}.
  The derivative of $\pi$ at $y$ is the
orthogonal projection from $T_yM$ onto $T_zM$. 
By the angle estimate above, this projection retains at least the
fraction $\alpha$ of the length of every tangent vector. Consequently,
its $m$-dimensional volume factor---the tangent Jacobian $J_m\pi(y)$---is
at least $\alpha^m$. Indeed, this Jacobian is the product of the $m$
singular values of the derivative, each of which is at least $\alpha$.

Since $\pi$ is injective, the area formula~\cite[Chapter~2, \S4]{Simon2014} gives
\[
\operatorname{Vol}(D)
=
\int_{M\cap B(z,\rho)}
J_mP(y)\,d\operatorname{Vol}_M(y)
\ge
\alpha^m\operatorname{Vol}\bigl(M\cap B(z,\rho)\bigr).
\]
Hence
\[
\operatorname{Vol}\bigl(M\cap B(z,\rho)\bigr)
\le
\alpha^{-m}\operatorname{Vol}(D)
\le
\frac{v_m\rho^m}
{\left(1-\frac{2\rho}{\tau}\right)^{m/2}}.\qedhere
\]
\end{proof}

\section{Curves in $\R^3$ and $\R^2$.}
\label{app:curves}
\begin{theorem}
  Let $M$ be a $1$-manifold without boundary in the unit cube $[0,1]^3$.
  For sufficiently many points, $\Eout<\Ein$ holds eventually
  whenever $w\leq 0.623$.
\end{theorem}
\begin{proof}
We proceed in the same way (and with the same notation)
as in Theorem~\ref{thm:special_case_non_empty_interval}, analyzing the quotient $\frac{A}{B}$. Since $v_1=2$, we get that
\[\frac{A}{B}=\frac{128}{47}\pi\Vol(M)(cr)^2.\]
And since $cr\leq\frac{5}{48}\tau$, we have
\[\frac{A}{B}\leq 0.0296\pi\Vol(M)\tau^2.\]
To bound the term $\pi\Vol(M)\tau^2$, note that the latter is simply
the volume of the $t$-thickened curve $M$. Since that curve fits
into the cube $[-\tau,1+\tau]^3$, we have
\[\pi\Vol(M)\tau^2\leq (1+2\tau)^3.\]
Finally, the reach $\tau$ of $M$ is at most $\sqrt{3}/2$ because
the reach of a curve in $\R^3$ is upper bounded
by half the diameter of the curve (i.e., the maximal distance
between two points on curve), which is at most $\sqrt{3}$,
the diameter of the cube; see~\cite[Thm 2.2]{akcmrw-estimating}.
Hence, we get $\pi\Vol(M)\tau^2\leq (1+\sqrt{3})^3$,
so
\[\frac{A}{B}\leq 0.0296(1+\sqrt{3})^3<0.604\]
and hence
\[\frac{1}{1+\frac{A}{B}}\geq 0.623.\qedhere\]
\end{proof}

\begin{theorem}
  Let $M$ be a $1$-manifold without boundary in the unit cube $[0,1]^2$.
  For sufficiently many points, $\Eout<\Ein$ holds eventually
  whenever $w\leq0.599$.
\end{theorem}
\begin{proof}
The proof is analogous to the previous one. 
Here, we got
\[\frac{A}{B}=\frac{48\pi}{47}\Vol(M)cr\leq \frac{5\pi}{47}\Vol(M)\tau.\]
$\Vol(M)\tau$ is half the area of the $\tau$-thickening of the curve $M$
which fits into the square $[-\tau,1+\tau]^2$, hence
\[2\Vol(M)\tau\leq(1+2\tau)^2\leq 4\]
where the last inequality comes from the fact that a closed
curve in the unit square
must enclose a tangential inscribed circle,
touching two points of $M$, and hence $\tau\leq\frac{1}{2}$. This yields
\[\frac{A}{B}\leq \frac{10\pi}{47}\]
and
\[\frac{1}{1+\frac{A}{B}}\geq 0.599.\qedhere\]
\end{proof}
\end{document}